\documentclass[preprint,12pt,3p,dvipsnames]{elsarticle}

\usepackage{amsmath, amssymb, amsfonts, amsbsy, amsthm, latexsym, epsfig, txfonts, palatino, pifont, tikz, color, graphicx, float}
\usetikzlibrary {positioning}
\definecolor {processblue}{cmyk}{0.96,0,0,0}

\theoremstyle{definition}

\newcommand{\longeq}{\scalebox{6}[1]{=}}

\journal{Bulletin of Mathematical Biosciences}

\makeatletter
\g@addto@macro\@floatboxreset\centering
\makeatother

\begin{document}

\begin{frontmatter}

\title{Local Immunodeficiency: Role of Neutral Viruses}

\author[label1]{Leonid Bunimovich}
\address[label1]{School of Mathematics,
Georgia Institute of Technology,
Atlanta, GA 30332-0160 USA}
\ead{bunimovh@math.gatech.edu}

\author[label2]{Longmei Shu\corref{cor1}}
\address[label2]{Mathematics \& Science Center, Emory University,
Suite W401,
400 Dowman Drive,
Atlanta, GA 30322}

\cortext[cor1]{Corresponding author}

\ead{lshu5@emory.edu}

\begin{abstract}
This paper analyzes the role of neutral viruses in the phenomenon of local immunodeficiency. We show that, even in the absence of altruistic viruses, neutral viruses can support the existence of persistent viruses, and thus local immunodeficiency. However, in all such cases neutral viruses can maintain only bounded (relatively small)  concentration of persistent viruses. Moreover, in all such cases the state of local immunodeficiency could only be marginally stable, while it is known that altruistic viruses can maintain stable local immunodeficiency. We also present an absolutely minimal cross-immunoreactivity network where a stable and robust state of local immunodeficiency can be maintained. It is now a challenge to synthetic biology to build such small networks with stable local immunodeficiency. Another important
{challenge for biology is to understand}
which types of viruses can play a role of persistent, altrustic and neutral ones, and whether a role which a given virus plays depends on the structure (topology) of a given cross-immunoreactivity network.
\end{abstract}


\end{frontmatter}




\section{Introduction}
Local immunodeficiency (LI) is a recently discovered phenomenon \cite{pnas} that appears in diseases characterised by cross-immunoreactivity of the corresponding pathogens (viruses). Examples of such diseases include Hepatitis C, HIV, dengue, influenza, etc
\cite{hattori,campo,nowak00,yoshioka,wodarz,nowak90,nowak91,nowak91s}. The phenomenon of local immunodeficiency means that some persistent antigens (viruses) manage to escape immune response because they are protected by altruistic viruses that take virtually all the response of the host's immune system on themselves. This discovery was made through the (numerical) analysis of a new model of Hepatitis C dynamics (evolution) that explained clinical and experimental observations that previous evolution models (and theory) of Hepatitis C failed to explain. Remarkably this new mathematical model 
{has}
fewer (types of) variables than the previous ones \cite{wodarz}. Yet this new model makes much more delicate exploration of the well-known phenomenon of cross-immunoreactivity than all previous evolution models of infectious diseases, including the fundamental dynamics model of HIV \cite{nowak90,nowak91,nowak91s}. Namely, 
{the new model}
does not assume that all cross-immunoreactivity interactions between different antigens (viruses) have the same (equal) strength. This was experimentally verified for Hepatitis C in
{the Center for Disease Control and Prevention (CDC)}
\cite{campo,hattori} before a final structure of a new model was created \cite{pnas}. 
{It is important to emphasize that all previous models of various diseases with cross-immunoreactivity always used the so called mean-field approximation where the strengths of interactions between all antigens (viruses) were assumed to be equal to one and the same constant.}

A striking discovery made in analyzing the dynamics of
{the new}
model was that
{in equilibrium state}
all intrahost viruses fall into one of three classes. The first class consists of persistent viruses that have the highest concentrations but the immune response against them is virtually zero. Therefore the  host's immune system
{"does not see" persistent viruses and thus}
demonstrates immunodeficiency against
{them}
. The second class of so called altruistic viruses is characterized by extremely low (virtually zero) concentrations.
{However}
almost all strength of the
{host's}
immune response goes
{against}
these altruistic viruses. This
{situation}
becomes possible because of special locations and structure of connections between persistent and altruistic viruses in the intrahost cross-immunoreactivity (CR) network of viruses (antigens) \cite{pnas}. Each of these two types of viruses comprise a very small (a few percent) part of all intrahost viruses. The rest of viruses (which comprise about 90\% of all) are called neutral.
{It is important to mention that such structure of cross-immunoreactivity networks appeared in all several hundred numerical experiments with our model \cite{pnas}. Therefore it seems to be an extremely stable and robust phenomenon.}

The host's immune system demonstrates immunodeficiency against persistent viruses because of their special positions (locations) within the intrahost CR network (CRN) and 
{because of the structure of}
their connections to altruistic viruses \cite{pnas}. This is where the term "local immunodeficiency" comes from.
{A central biological question in this area is to determine which types of viruses can play a role of persistent and/or altruistic viruses. Indeed, it's very interesting biologically and important for public health policies on how to handle the diseases with cross-immunoreactivity}

Paper \cite{BS} shows that the phenomenon of local immunodeficiency is typically stable and robust under various realistic conditions. Moreover, it 
{was proved in this paper}
that stable and robust local immunodeficiency can already occur in very small CR networks consisting of just three nodes.
{We already know diseases like Hepatitis C have large CR networks, now we have shown small CR networks can also have stable local immunodeficiency. These results are very interesting in general and have attracted experimentalists. For example, the recent paper \cite{NM} studied a social structure of population of viruses in an experimental framework. Interestingly, their results showed a different type of altruism.}
Indeed,
{true}
altruists should sacrifice themselves for others (as \cite{pnas} demonstrated) but the viruses called altruistic in \cite{NM} help themselves and also some other viruses. It is not
{quite}
altruism
{as it is usually understood.}

In the present paper we analyze the role of neutral viruses in CR networks and their ability to maintain persistent viruses, i.e. to generate local immunodeficiency. The question whether and why neutral viruses are needed to create and maintain local immunodeficiency remained unanswered in previous studies. We show that an answer to this natural question is nontrivial and unexpected (as is essentially everything in the studies of the phenomenon of local immunodeficiency so far). Namely
{we show that}
neutral viruses (without the presence of altruistic viruses) could maintain only marginally stable state of a local immunodeficiency. Moreover, without altruistic viruses the population of persistent viruses can not be large and is bounded from above by exact computed values.
{This contrasts with the results of \cite{BS} that in the presence of altruistic viruses there is no upper bound on the population size of persistent viruses.}

{Another fundamental question is whether some viruses may play different roles becoming persistent, altruistic or neutral depending on the structure of a smaller subnetwork that belongs to a (large) cross-immunoreactivity network. Also altruistic viruses should be a primary target of public health policies.}

We also present here a new minimal (three-node) CRN with a stable and robust state of a local immunodeficiency. This
{new}
one and the one found in \cite{BS} are the only two CRNs with just three types of viruses which maintain a stable and robust local imunnodeficiency.
{However, the three viruses network presented here is in fact the absolutely minimal because it contains just two edges, the minimal number of edges in any network with three nodes.}
It would be very interesting to create small cross-immunoreactivity networks which demonstrate stable and robust phenomenon of local immunodeficiency. This seems to be a natural challenge to synthetic biology,
{and the absolutely minimal network presented in this paper is the first candidate to be considered in such studies because it has the the minimal number of elements (nodes) and the simplest possible structure (topology) of a network.}

\section{Model of evolution of a disease with heterogeneous CR network}\label{model}
In this section we define the model of the 
{Hepatitis C (HC)}
evolution introduced in \cite{pnas}. It is important to 
{note}
that this model is applicable to any disease with cross-immunoreactivity.

Consider a 
{system}
which contains a population of $n$ viral antigenic variants $x_i$ inducing $n$ immune responses $r_i$ in the form of antibodies (Abs). The viral variants exhibit cross-immunoreactivity which results in 
a CR network.  The latter can be represented as a directed graph $G_{CRN}=(V,E)$, with vertices corresponding to viral variants and directed edges connecting CR variants.
Because not all interactions with Ab lead to neutralization, we consider two sets of weight functions for the CRN. These functions are defined by immune neutralization and immune stimulation matrices $U=(u_{ij})_{i,j=1}^n$ and $V=(v_{ij})_{i,j=1}^n$, where $0\le u_{ij},v_{ij}\le1$; $u_{ij}$ represents the binding affinity of Ab to $j~(r_j)$ with the $i$-th variant; and $v_{ij}$ reflects the strength of stimulation of Ab to $j~(r_j)$ by the $i$-th variant. The immune response $r_i$ against variant $x_i$ is neutralizing; i.e., $u_{ii}=v_{ii}=1$. The 
{resulting}
evolution of the antigens (viruses) and antibodies populations is given by the following system of ordinary differential equations \cite{pnas,BS}.

\begin{equation}\label{population}
\begin{split}
\dot x_i=f_ix_i-px_i\sum_{j=1}^nu_{ji}r_j,\quad i=1,\dots,n,\\
\dot r_i=c\sum_{j=1}^nx_j\frac{v_{ji}r_i}{\sum_{k=1}^nv_{jk}r_k}-br_i,\quad i=1,\dots,n.
\end{split}
\end{equation}

{In this model}
the viral variant $x_i$ replicates at the rate $f_i$ and is eliminated by the immune responses $r_j$ at the rates $pu_{ji}r_j$,
{where $p$ is a constant}.
The immune responses $r_i$ are stimulated by the $j$-th variant at the rates $cg_{ji}x_j$, where $g_{ji}=\frac{v_{ji}r_i}{\sum_{k=1}^nv_{jk}r_k}$ represents the probability of stimulation of the immune response $r_i$ by the variant $x_j$ 
{and $c$ is a constant}.
This model \cite{pnas} allows us to incorporate the phenomenon of the original antigenic sin  \cite{tom,pan,shin,hvi,jin,alotofauthors}, which states that $x_i$ preferentially stimulates preexisting immune responses capable of binding to $x_i$. The immune response $r_i$ decays at rate $b$ in the absence of stimulation.

{As in \cite{pnas,BS}}
we consider the situation where the immune stimulation and neutralization coefficients are equal to constants $\alpha$ and $\beta$, respectively. To be more specific, both the immune neutralization and stimulation matrices are completely defined by the structure of the CR network, i.e.,
$$U=\text{Id}+\beta A^T,V=\text{Id}+\alpha A,$$
where $A$ is the adjacency matrix of $G_{CRN}$. In the absence of CR among viral variants the system  reduces to the model developed in \cite{nowak00} for a heterogeneous viral population.
Because the neutralization of an antigen may require more than one antibodies, we assume that $0<\beta=\alpha^k<\alpha<1$ \cite{pnas}. It is important to mention that we analyze a more general model here than the one studied in \cite{pnas}, where it was assumed that all viruses replicate at the same rate
{, i.e. all the $f_i$'s are equal to each other. Here, as in \cite{BS} we allow the $f_i$'s to be different and consider all possible values of them as long as they make sense, i.e. are positive.}

\section{An absolutely minimal network with stable local immunodeficiency}

{In this section we present a minimal cross-immunoreactivity network with stable and robust state of local immunodeficiency.}

{Recall that an equilibrium solution to a system of differential equations is called \textbf{stable} if the solution to the system under any perturbation of the \textbf{initial condition} converges with time to this solution. A solution is \textbf{robust} if there is a solution to the system under small perturbations of its \textbf{parameters} (coefficients) with dynamics  analogous to the solution of the unperturbed system.}

{As in \cite{BS} we will consider "extreme" stable fixed points (equilibrium solutions) with local immunodeficiency. This means that we assume the immune response against a persistent virus is zero and the concentration of an altruistic virus is also identically zero. Of course, these may sound like unrealistic assumptions. Observe, however, that the right hand side of system \eqref{population} depends on parameters continuously (and even smoothly); under small perturbations of parameters there will be a stable fixed point close to the extreme one with stable local immunodeficiency (LI) in the unperturbed system of ODEs (see e.g. \cite{Chicone}).}

{Certainly these stable solutions should make (biological) sense, i.e. both the  concentrations of altruistic viruses and the immune response against persistent viruses have to be positive. Both these values will be close to zero under small perturbations of parameters of system \eqref{population}. Fixed points that make sense do exist because a neighborhood of zero (here the fixed points depend continuously on parameters of the system) always contains a subset with only positive values for all coordinates.}

{Observe also that in reality a system is never exactly at a fixed point but close to it because of ever-present small perturbations/fluctuations. The dynamics always pushes the sate of the system closer to the steady state. Therefore a real system will evolve towards a fixed point. This evolution takes infinite time and the steady state will never be achieved although (in absence of random perturbations) the current state of the system will become closer and closer to this stable fixed point with time.}

{It is  worthwhile to mention here that in numerical simulations one must set a (positive) threshold such that if the population size of some virus becomes less than this threshold then this virus is considered to be eliminated. For instance, simulations of system \eqref{population} in \cite{pnas} assumed that a virus gets eliminated by the immune system if its concentration becomes lower than the initial concentration of this virus.}

Consider the following network with 3 types of viruses.
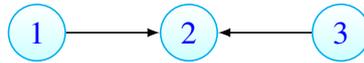
\begin {figure}[h!]
\centering
\begin {tikzpicture}[-latex ,auto ,node distance =02cm and 2cm ,on grid ,
semithick ,
state/.style ={ circle ,top color =white , bottom color = processblue!20 ,
draw,processblue , text=blue , minimum width =0.5 cm}]
\node[state] (A)
{$1$};
\node[state] (B) [right=of A] {$2$};
\node[state] (C) [right =of B] {$3$};
\path (A) edge [] node[below =0.15 cm] {} (B);
\path (C) edge [] node[below =0.15 cm] {} (B);
\end{tikzpicture}
\caption{Symmetric CRN}
\label{symmetric}
\end{figure}

{Although this network looks rather similar to the three viruses network with stable and robust local immunodeficiency found in \cite{BS}, it actually has fewer edges. We will show that this smaller cross-immunoreactivity network is also able to support stable and robust local immunodeficiency.}

The adjacency matrix and the corresponding neutralization and stimulation matrices
{of this network}
are
$$A=\begin{pmatrix}0 & 1 & 0 \\ 0 & 0 & 0 \\ 0 & 1 & 0\end{pmatrix}, U=\begin{pmatrix}1 & 0 & 0 \\ \beta & 1 & \beta \\ 0 & 0 & 1\end{pmatrix}, V=\begin{pmatrix}1 & \alpha & 0 \\ 0 & 1 & 0 \\ 0 & \alpha & 1\end{pmatrix}.$$
{All our networks are finite simple graphs, therefore the adjacency matrix has 0-1 entries indicating whether there is an edge going from node i to node j.}

The model of population evolution for this network is
{given by the following differential equations}
$$\begin{cases}
\dot x_1=f_1x_1-px_1(r_1+\beta r_2),\\
\dot x_2=f_2x_2-px_2r_2,\\
\dot x_3=f_3x_3-px_3(\beta r_2+r_3),\\
\dot r_1=c\frac{x_1r_1}{r_1+\alpha r_2}-br_1,\\
\dot r_2=c(\frac{\alpha x_1r_2}{r_1+\alpha r_2}+\frac{x_2r_2}{r_2}+\frac{\alpha x_3r_2}{\alpha r_2+r_3})-br_2,\\
\dot r_3=c\frac{x_3r_3}{\alpha r_2+r_3}-br_3.
\end{cases}$$
The Jacobian of this model is
$$J=\begin{pmatrix}A_J & B \\ C & D\end{pmatrix},
{\text{ where }}
A_J=\begin{pmatrix}f_1-p(r_1+\beta r_2) & 0 & 0 \\ 0 & f_2-pr_2 & 0 \\ 0 & 0 & f_3-p(\beta r_2 + r_3)\end{pmatrix},$$
$$B=\begin{pmatrix}-px_1 & -p\beta x_1 & 0 \\ 0 & -px_2 & 0 \\ 0 & -p\beta x_3 & -px_3\end{pmatrix}, C=\begin{pmatrix} \frac{cr_1}{r_1+\alpha r_2} & 0 & 0 \\ \frac{c\alpha r_2}{r_1+\alpha r_2} & c & \frac{c\alpha r_2}{\alpha r_2+r_3} \\ 0 & 0 & \frac{cr_3}{\alpha r_2+r_3}\end{pmatrix},$$
$$D=\begin{pmatrix}\frac{cx_1\alpha r_2}{(r_1+\alpha r_2)^2}-b & -\frac{c\alpha x_1r_1}{(r_1+\alpha r_2)^2} & 0 \\ -\frac{c\alpha x_1r_2}{(r_1+\alpha r_2)^2} & \frac{c\alpha x_1r_1}{(r_1+\alpha r_2)^2}+\frac{c\alpha x_3r_3}{(\alpha r_2+r_3)^2}-b & -\frac{c\alpha x_3r_2}{(\alpha r_2+r_3)^2} \\ 0 & -\frac{c\alpha x_2r_3}{(\alpha r_2+r_3)^2} & \frac{cx_3\alpha r_2}{(\alpha r_2+r_3)^2}-b\end{pmatrix}.$$
Based on the observation that altruistic nodes should have the highest in-degree \cite{pnas,BS}, we will be looking for a fixed point 
{with local immunodeficiency}
where node 2 is altruistic 
{(no virus concentration but high immune response, i.e. $x_i=0,r_i>0$).}
Since the nodes 1 and 3 are symmetric, without loss of generality we pick
{either}
one of them to be persistent 
{(high virus concentration but no immune response, i.e. $x_i>0,r_i=0$)}
and the other one to be neutral 
{(positive virus concentration and immune response, i.e. $x_i>0,r_i>0$)}
, i.e.
$$x_2=0,r_2>0,x_1>0,r_1=0,x_3>0,r_3>0.$$

{The virus populations and corresponding antibody populations at such a fixed point (steady state) with local immunodeficiency are}
$$r_1=0,r_2=\frac{f_1}{p\beta},r_3=\frac{f_3-f_1}{p}>0,f_3>f_1.$$
$$x_1=(br_2-\frac{c\alpha x_3r_2}{\alpha r_2+r_3})\frac{r_1+\alpha r_2}{c\alpha r_2}=\frac{b}{c}r_2(1-\alpha),x_2=0,x_3=\frac{b}{c}(\alpha r_2+r_3).$$

{To see whether this fixed point with local immunodeficiency we just found is stable, we consider the Jacobian of the system at this point.}

$$A_J=\begin{pmatrix}0 & 0 & 0 \\ 0 & f_2-pr_2 & 0 \\ 0 & 0 & 0\end{pmatrix}, B=\begin{pmatrix}-px_1 & -p\beta x_1 & 0 \\ 0 & 0 & 0 \\ 0 -p\beta x_3 & -px_3\end{pmatrix},$$
$$C=\begin{pmatrix}0 & 0 & 0 \\ c & c & \frac{c\alpha r_2}{\alpha r_2+r_3} \\ 0 & 0 & \frac{cr_3}{\alpha r_2+r_3}\end{pmatrix}, D=\begin{pmatrix}\frac{cx_1}{\alpha r_2}-b & 0 & 0 \\ -\frac{cx_1}{\alpha r_2} & \frac{c\alpha x_3r_3}{(\alpha r_2+r_3)^2}-b & -\frac{c\alpha x_3r_2}{(\alpha r_2+r_3)^2} \\ 0 & -\frac{c\alpha x_3r_3}{(\alpha r_2+r_3)^2} & \frac{cx_3\alpha r_2}{(\alpha r_2+r_3)^2}-b\end{pmatrix},$$
$$J=\begin{pmatrix}0 & 0 & 0 & -px_1 & -p\beta x_1 & 0 \\ 0 & f_2-pr_2 & 0 & 0 & 0 & 0 \\ 0 & 0 & 0 & 0 & -p\beta x_3 & -px_3 \\ 0 & 0 & 0 & \frac{b}{\alpha}-2b & 0 & 0 \\ c & c & \frac{c\alpha r_2}{\alpha r_2+r_3} & b-\frac{b}{\alpha} & \frac{b\alpha r_3}{\alpha r_2+r_3}-b & -\frac{b\alpha r_2}{\alpha r_2+r_3} \\ 0 & 0 & \frac{cr_3}{\alpha r_2+r_3} & 0 & -\frac{b\alpha r_3}{\alpha r_2+r_3} & -\frac{br_3}{\alpha r_2+r_3}\end{pmatrix}.$$

{If all the eigenvalues of the Jacobian are negative, or complex with negative real parts, then the fixed point has stable local immunodeficiency. Since the eigenvalues are the roots of the characteristic polynomial, we need to compute the characteristic polynomial of the Jacobian matrix.}

Let $\lambda_1=f_2-pr_2=f_2-f_1/\beta,\lambda_2=\frac{b}{\alpha}-2b$.
$$\det(J-\lambda I)=(\lambda-\lambda_1)(\lambda-\lambda_2)\begin{vmatrix} -\lambda & 0 & -p\beta x_1 & 0 \\ 0 & -\lambda & -p\beta x_3 & -px_3 \\ c & \frac{c\alpha r_2}{\alpha r_2+r_3} & \frac{b\alpha r_3}{\alpha r_2+r_3}-b-\lambda & -\frac{b\alpha r_2}{\alpha r_2+r_3} \\ 0 & \frac{cr_3}{\alpha r_2+r_3} & -\frac{b\alpha r_3}{\alpha r_2+r_3} & -\frac{br_3}{\alpha r_2+r_3}-\lambda \end{vmatrix}$$
$$=(\lambda-\lambda_1)(\lambda-\lambda_2)P(\lambda).$$
{Detailed computation of $P(\lambda)$ is shown in the appendix. The result is}
$$P(\lambda)=\alpha bf_1\lambda^2+bf_1(1-\alpha)(\lambda^2+\frac{br_3}{\alpha r_2+r_3}\lambda+bpr_3)+\lambda^2(\lambda+b)(\lambda+\frac{br_3(1-\alpha)}{\alpha r_2+r_3})+px_3\lambda[\frac{cr_3}{\alpha r_2+r_3}\lambda+\frac{bcr_3(1-\alpha)}{\alpha r_2+r_3}].$$
{Observe that}
all coefficients of the polynomial $P(\lambda)$ are positive. Therefore it cannot have positive real roots.
Combining all conditions on the parameters  $f_3>f_1,\lambda_1=f_2-f_1/\beta<0,\lambda_2=b/\alpha-2b<0,$ one gets
$\beta f_2<f_1<f_3,\alpha>1/2$. 

We present now some numerical examples where the LI is stable.
{These numerical values are consistent with the ranges of the corresponding parameters considered in the literature \cite{nowak00, rong, dahari, pnas, BS}.}
\begin{enumerate}
    \item $f_1=2,f_2=3,f_3=3,p=2,c=1,b=3,\alpha=2/3,\beta=4/9$. With these parameters, $\lambda_1=\lambda_2=-1.5$, the polynomial  $P(\lambda)$ has 2 pairs of complex roots, both with negative real parts.
   
    \item $f_1=2,f_2=2,f_3=3,p=2,c=1,b=1,\alpha=3/4,\beta=9/16$. Under such parameters, $\lambda_1=-14/9<0,\lambda_2=-2/3<0$, $P(\lambda)$ has 2 pairs of complex roots, both with negative real parts.
\end{enumerate}

By continuity this fixed point is stable on a positive measure set in the parameter space.
{This means that this fixed point is physically observable, i.e. there is a positive probability (real life chance) to generate (e.g. numerically) a system with the CRN in Fig. \ref{symmetric} which has a stable and robust steady state of local immunodeficiency. It (positive measure in parameter space) also means that there is a real chance to generate a biological CR viral network with stable local immunodeficiency.}

{Observe that we only computed one fixed point. Certainly there are other stable fixed points with local immunodeficiency that are stable and robust under perturbation of initial conditions and parameters. For example, there is at least one other fixed point where we switch the roles of node 1 and 3. Also small variations of parameters result in a "shift" of a fixed point. Such fixed points will remain stable in view of continuity of the coefficients of the system of differential equations under study.}

\begin {figure}[h]
\begin {tikzpicture}[-latex ,auto ,node distance =2 cm and 2cm ,on grid ,
semithick ,
state/.style ={ circle ,top color =white , bottom color = processblue!20 ,
draw,processblue , text=blue , minimum width =0.5 cm}]
\node[state] (A)
{$1$};
\node[state] (B) [right=of A] {$2$};
\node[state] (C) [right =of B] {$3$};
\path (A) edge [] node[below =0.15 cm] {} (B);
\path (B) edge [bend left =15] node[above] {} (C);
\path (C) edge [bend left =15] node[below] {} (B);
\end{tikzpicture}
\caption{branch-cycle CRN}
\label{minimal}
\end{figure}
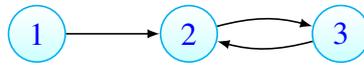

Compare this network in Fig. \ref{symmetric} to the network with stable and robust local immunodeficiency found in \cite{BS} (Fig. \ref{minimal}).

{It is easy to see that the new network in Fig. \ref{symmetric} is really simpler than the branch-cycle CRN.}
Indeed both networks contain three types of viruses but the symmetric CRN contains fewer edges.
{This shows that the removal of an edge going out of an altruistic node (virus) may not destroy the state of stable local immunodeficiency.}

\section{A role of neutral viruses in CR networks}

{In this section we will analyze}
the role of neutral viruses in local immunodeficiency, i.e. whether neutral viruses alone, without altruistic viruses, can sustain stable local immunodeficiency.
{Recall}
that persistent nodes represent types of viruses whose concentration is high but immune response against them is zero, and neutral nodes are the ones where both the virus and the antibody population are positive. From now on we 
{will}
use the
{light pink}
color to represent persistent nodes and the
{dark}
green color for neutral nodes in
{CR networks}.

Consider first the simplest CR network consisting of just one persistent and one neutral node. It is an asymmetric network, where the persistent node is connected to the neutral one.
Dynamics equations for evolution of this system are

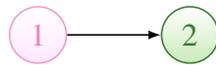
\begin {figure}[h]
\begin {tikzpicture}[-latex ,auto ,node distance =2 cm and 2cm ,on grid ,
semithick ,
per/.style={circle,top color=white,bottom color=Lavender!20,draw,Lavender,text=CarnationPink,minimum width=0.5cm},
alt/.style={circle,top color=white,bottom color=processblue!20,draw,processblue,text=blue,minimum width=0.5cm},
nac/.style={circle,top color=white, bottom color=OliveGreen!20,draw,OliveGreen,text=OliveGreen,minimum width=0.5cm},
idl/.style ={circle,top color=white,bottom color=gray!20,draw,gray,text=black,minimum width =0.5 cm}]
\node[per] (A) {$1$};
\node[nac] (B) [right =of A] {$2$};
\path (A) edge [] node[above] {} (B);
\end{tikzpicture}
\caption{size 2 CRN}
\label{2fig}
\end{figure}

$$\begin{cases}
\dot x_1=f_1x_1-px_1(r_1+\beta r_2)\\
\dot x_2=f_2x_2-px_2r_2\\
\dot r_1=cx_1\frac{r_1}{r_1+\alpha r_2}-br_1\\
\dot r_2=c(x_1\frac{\alpha r_2}{r_1+\alpha r_2}+x_2)-br_2
\end{cases}$$
A family of fixed points where node 1 is persistent and node 2 is neutral is given by the following relations.
$$f_1=\beta f_2,0<x_1<\frac{bf_2}{cp},x_2=\frac{bf_2}{cp}-x_1,r_1=0,r_2=\frac{f_2}{p}$$
{Notice that this represents infinitely many numbers of fixed points (equilibrium states) that vary continuously on a line segment. To study the stability of these fixed points we need to find the eigenvalues of the Jacobian.}
The Jacobian of the system is 
$$J=\begin{pmatrix}
f_1-p(r_1+\beta r_2) & 0 & -px_1 & -p\beta x_1\\
0 & f_2-pr_2 & 0 & -px_2\\
\frac{cr_1}{r_1+\alpha r_2} & 0 & \frac{cx_1\alpha r_2}{(r_1+\alpha r_2)^2}-b & -\frac{cx_1\alpha r_1}{(r_1+\alpha r_2)^2}\\
\frac{c\alpha r_2}{r_1+\alpha r_2} & c & -\frac{cx_1\alpha r_2}{(r_1+\alpha r_2)^2} & \frac{cx_1\alpha r_1}{(r_1+\alpha r_2)^2}-b
\end{pmatrix},$$

At those fixed points the Jacobian becomes 
$$J=\begin{pmatrix}0 & 0 & -px_1 & -p\beta x_1\\
0 & 0 & 0 & -px_2\\
0 & 0 & \frac{cx_1}{\alpha r_2}-b & 0\\
c & c & -\frac{cx_1}{\alpha r_2} & -b\end{pmatrix}.$$ Let $\lambda_1=\frac{cx_1}{\alpha r_2}-b$.
{The detailed computation of the characteristic polynomial of the Jacobian is shown in the appendix, the result is}
$$\det(J-\lambda I)=\lambda(\lambda-\lambda_1)[\lambda^2+b\lambda+cp(x_2+\beta x_1)]=\lambda(\lambda-\lambda_1)P(\lambda),$$

where
$$\lambda_1=\frac{cx_1}{\alpha r_2}-b=\frac{cp}{\alpha f_2}x_1-b<\frac{b}{\alpha}-b.$$
All coefficients of the quadratic polynomial $P(\lambda)$ are positive. Therefore its roots are either real, negative or complex with negative real parts. Depending on the value of $x_1$ the eigenvalue $\lambda_1$ could be positive or negative. And for every fixed point in the family, the Jacobian has a 0 eigenvalue. This means that the corresponding state (fixed point with local immunodeficiency) is never stable but it could be marginally stable (because of the zero eigenvalue). However in this case the concentration of persistent viruses cannot exceed some fixed value. 

Indeed for this family of fixed points, when $x_1$ is small ($x_1<\alpha\frac{bf_2}{cp}$) the fixed points are stable on the subspace $f_1=\beta f_2$; and when $x_1$ is big ($\alpha\frac{bf_2}{cp}<x_1<\frac{bf_2}{cp}$) the fixed points are unstable. It should be contrasted with the results of \cite{BS} where, in the presence of an altruistic node, the state of local immunodeficiency can be stable (rather than marginally stable) and there is no bound on the concentration of persistent viruses.

Consider now a larger network with one persistent and two neutral nodes.

\begin{figure}[H]
\begin{tikzpicture}[-latex ,auto ,node distance =02cm and 2cm ,on grid ,
semithick ,
per/.style={circle,top color=white,bottom color=Lavender!20,draw,Lavender,text=CarnationPink,minimum width=0.5cm},
alt/.style={circle,top color=white,bottom color=processblue!20,draw,processblue,text=blue,minimum width=0.5cm},
nac/.style={circle,top color=white, bottom color=OliveGreen!20,draw,OliveGreen,text=OliveGreen,minimum width=0.5cm},
idl/.style ={circle,top color=white,bottom color=gray!20,draw,gray,text=black,minimum width =0.5 cm}]
\node[per] (A) {$1$};
\node[nac] (B) [right =of A] {$2$};
\node[nac] (C) [left =of A] {$3$};
\path (A) edge [] node[above] {} (B);
\path (A) edge [] node[above] {} (C);
\end{tikzpicture}
\label{3-node}
\caption{size 3 CRN}
\end{figure}
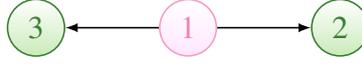

The adjacency matrix, neutralization and stimulation matrices
{for this CR network}
are
$$A=\begin{pmatrix}
0 & 1 & 1\\
0 & 0 & 0\\
0 & 0 & 0
\end{pmatrix}, 
U=\begin{pmatrix}
1 & 0 & 0\\
\beta & 1 & 0\\
\beta & 0 & 1
\end{pmatrix},
V=\begin{pmatrix}
1 & \alpha & \alpha \\
0 & 1 & 0\\
0 & 0 & 1
\end{pmatrix}.$$

The population growth equations for this system are
$$\begin{cases}
\dot x_1=f_1x_1-px_1(r_1+\beta r_2+\beta r_3)\\
\dot x_2=f_2x_2-px_2r_2\\
\dot x_3=f_3x_3-px_3r_3\\
\dot r_1=cx_1\frac{r_1}{r_1+\alpha r_2+\alpha r_3}-br_1\\
\dot r_2=cx_1\frac{\alpha r_2}{r_1+\alpha r_2+\alpha r_3}+cx_2-br_2\\
\dot r_3=cx_1\frac{\alpha r_3}{r_1+\alpha r_2+\alpha r_3}+cx_3-br_3
\end{cases}$$
A family of fixed points for this network is given via the following relations 
$$f_1=\beta(f_2+f_3),0<x_1<\frac{b}{c}(r_2+r_3),x_2=\frac{b}{c}r_2-x_1\frac{r_2}{r_2+r_3},x_3=\frac{b}{c}r_3-x_1\frac{r_3}{r_2+r_3},$$
$$r_1=0,r_2=\frac{f_2}{p},r_3=\frac{f_3}{p}.$$
{Again this is a group of steady states with local immunodeficiency (node 1 is persistent with positive virus population but 0 antibody population) that vary continuously on a line segment.}
The Jacobian of the system is
$$J=\begin{pmatrix}f_1-p(r_1+\beta r_2+\beta r_3) & 0 & 0 & -px_1 & -p\beta x_1 & -p\beta x_1\\
0 & f_2-pr_2 & 0 & 0 & -px_2 & 0 \\ 0 & 0 & f_3-pr_3 & 0 & 0 & -px_3 \\
\frac{cr_1}{r_1+\alpha r_2+\alpha r_3} & 0 & 0 & \frac{cx_1\alpha(r_2+r_3)}{(r_1+\alpha r_2+\alpha r_3)^2}-b & -\frac{cx_1\alpha r_1}{(r_1+\alpha r_2+\alpha r_3)^2} & -\frac{cx_1\alpha r_1}{(r_1+\alpha r_2+\alpha r_3)^2} \\
\frac{c\alpha r_2}{r_1+\alpha r_2+\alpha r_3} & c & 0 & -\frac{cx_1\alpha r_2}{(r_1+\alpha r_2+\alpha r_3)^2} & \frac{cx_1\alpha(r_1+\alpha r_3)}{(r_1+\alpha r_2+\alpha r_3)^2}-b & -\frac{cx_1\alpha^2r_2}{(r_1+\alpha r_2+\alpha r_3)^2} \\
\frac{c\alpha r_3}{r_1+\alpha r_2+\alpha r_3} & 0 & c & -\frac{cx_1\alpha r_3}{(r_1+\alpha r_2+\alpha r_3)^2} & -\frac{cx_1\alpha^2r_3}{(r_1+\alpha r_2+\alpha r_3)^2} & \frac{cx_1\alpha(r_1+\alpha r_2)}{(r_1+\alpha r_2+\alpha r_3)^2}-b \end{pmatrix}.$$
At the fixed points under analysis we have
$$J=\begin{pmatrix} 0 & 0 & 0 & -px_1 & -p\beta x_1 & -p\beta x_1 \\
0 & 0 & 0 & 0 & -px_2 & 0 \\
0 & 0 & 0 & 0 & 0 & -px_3 \\
0 & 0 & 0 & \frac{cx_1}{\alpha(r_2+r_3)}-b & 0 & 0 \\
\frac{cr_2}{r_2+r_3} & c & 0 & -\frac{cx_1r_2}{\alpha(r_2+r_3)^2} & \frac{cx_1r_3}{(r_2+r_3)^2}-b & -\frac{cx_1r_2}{(r_2+r_3)^2} \\
\frac{cr_3}{r_2+r_3} & 0 & c & -\frac{cx_1r_3}{\alpha(r_2+r_3)^2} & -\frac{cx_1r_3}{(r_2+r_3)^2} & \frac{cx_1r_2}{(r_2+r_3)^2}-b\end{pmatrix}.$$

{We compute the characteristic polynomial of the Jacobian for stability analysis of the fixed points.
Let $\lambda_1=\frac{cx_1}{\alpha(r_2+r_3)}-b<\frac{b}{\alpha}-b$, detailed computation is put in the appendix and the characteristic polynomial is}
$$|J-\lambda I|=(\lambda-\lambda_1)[\lambda D_1(\lambda)+px_2 D_2(\lambda)],$$
where
$$D_1(\lambda)=\lambda(\lambda+b-\frac{cx_1}{r_2+r_3})[\lambda^2+b\lambda+cp\beta x_1]+cpx_3[\lambda^2+(b-\frac{cx_1r_3}{(r_2+r_3)^2})\lambda+cp\beta x_1\frac{r_2}{r_2+r_3}],$$
$$D_2(\lambda)=c\lambda[\lambda^2+(b-\frac{cx_1r_2}{(r_2+r_3)^2})\lambda+cpx_3]+cp\beta x_1\lambda\frac{cr_3}{r_2+r_3}.$$
{So, in this case}
again 0 is an eigenvalue. Another eigenvalue $\lambda_1$ can be positive or negative depending on the value of $x_1$. The rest of the eigenvalues are the roots of a degree four polynomial whose coefficients are all positive.

Therefore for small values of $x_1$, namely if $x_1<\alpha\frac{b}{c}(r_2+r_3)$, the fixed points are stable on positive measure subsets of the subspace $f_1=\beta(f_2+f_3)$. If $x_1$ has larger values   $\alpha\frac{b}{c}(r_2+r_3)<x_1<\frac{b}{c}(r_2+r_3)$ then the fixed points are unstable.

In general, for any size network to have fixed points with only persistent and neutral nodes (viruses), the parameters $f_i$'s have to satisfy some condition that forms a positive codimension subspace \cite{BS}. And based on the previous two examples, it is natural to assume that there is a family of fixed points, where the Jacobian has a 0 eigenvalue, and an eigenvalue whose sign depends on the size of the viral population of the persistent node. The rest of the eigenvalues are the roots of a polynomial with positive coefficients.

\begin{figure}[h]
\begin{tikzpicture}[-latex ,auto ,node distance =02cm and 2cm ,on grid ,
semithick ,
per/.style={circle,top color=white,bottom color=Lavender!20,draw,Lavender,text=CarnationPink,minimum width=0.5cm},
alt/.style={circle,top color=white,bottom color=processblue!20,draw,processblue,text=blue,minimum width=0.5cm},
nac/.style={circle,top color=white, bottom color=OliveGreen!20,draw,OliveGreen,text=OliveGreen,minimum width=0.5cm},
idl/.style ={circle,top color=white,bottom color=gray!20,draw,gray,text=black,minimum width =0.5 cm}]
\node[per] (A) at (0,0) {$1$};
\node[nac] (B) at (2,0) {$2$};
\node[nac] (C) at (1.5,1.5) {$3$};
\node[nac] (D) at (0,2) {$4$};
\node[nac] (E) at (-1.5,1.5) {$5$};
\node[nac] (F) at (-2,0) {$6$};
\node[nac] (G) at (-1.5,-1.5) {$7$};
\node[nac] (H) at (1.5,-1.5) {$n$};
\node[text width=1pt] at (0,-2) {$\cdot$};
\node[text width=1pt] at (-0.4,-1.9) {$\cdot$};
\node[text width=1pt] at (0.4,-1.9) {$\cdot$};
\path (A) edge [] node[above] {} (B);
\path (A) edge [] node[above] {} (C);
\path (A) edge [] node [above] {} (D);
\path (A) edge [] node [above] {} (E);
\path (A) edge [] node [above] {} (F);
\path (A) edge [] node [above] {} (G);
\path (A) edge [] node [above] {} (H);
\end{tikzpicture}
\caption{size $n$ CRN}
\label{n-node}
\end{figure}
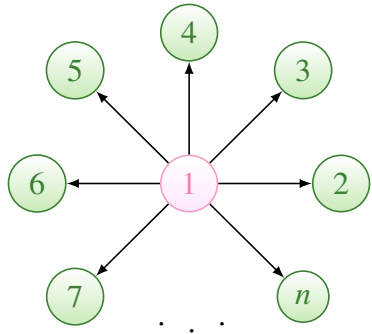

Consider now a CR network with one persistent node and any (finite) number  of neutral nodes to support it. Without loss of generality we assume that the persistent node is node 1, and the nodes 2 through $n$ are neutral (Figure \ref{n-node}). There is an edge going from the persistent node to each of the neutral nodes.

The adjacency matrix of such a network is
$$A=\begin{pmatrix}0 & 1 & \cdots & 1 \\
0 & 0 & \cdots & 0 \\
\vdots & \vdots & \ddots & \vdots \\
0 & 0 & \cdots & 0
\end{pmatrix}.$$
Then  $$U=\text{Id}+\beta A^T=\begin{pmatrix}
1 & 0 & \cdots & 0 \\
\beta & 1 & \cdots & 0 \\
\vdots & \vdots & \ddots & \vdots \\
\beta & 0 & \cdots & 1
\end{pmatrix}, V=\text{Id}+\alpha A=\begin{pmatrix}
1 & \alpha & \cdots & \alpha \\
0 & 1 & \cdots & 0 \\
\vdots & \vdots & \ddots & \vdots \\
0 & 0 & \cdots & 1
\end{pmatrix}.$$
The evolution equations for this system are
\vspace{1pt}\\
$$\begin{cases}
\dot x_1=f_1 x_1-px_1(r_1+\beta\sum_{j=2}^nr_j),\\
\vspace{1pt}\\
\dot x_i=f_ix_i-px_ir_i,2\le i\le n;\\
\vspace{1pt}\\
\dot r_1=c\dfrac{x_1r_1}{r_1+\alpha\sum_{j=2}^nr_j}-br_1,\\
\vspace{1pt}\\
\dot r_i=c\dfrac{\alpha x_1r_i}{r_1+\alpha\sum_{j=2}^nr_j}+cx_i-br_i,2\le i\le n.
\end{cases}$$
\vspace{1pt}\\

Let $N_r=\sum_{j=2}^nr_j$.   Consider the Jacobian of this system.

$$J=\begin{pmatrix}A_J & B \\ C & D\end{pmatrix},
A_J=\begin{pmatrix}
f_1-p(r_1+\beta N_r) & 0 & \cdots & 0 \\
0 & f_2-pr_2 & \cdots & 0 \\
\vdots & \vdots & \ddots & \vdots \\
0 & 0 & \cdots & f_n-pr_n
\end{pmatrix},$$
$$B=\begin{pmatrix}
-px_1 & -p\beta x_1 & \cdots & -p\beta x_1 \\
0 & -px_2 & \cdots & 0 \\
\vdots & \vdots & \ddots & \vdots \\
0 & 0 & \cdots & -px_n
\end{pmatrix},C=\begin{pmatrix}
c\frac{r_1}{r_1+\alpha N_r} & 0 & \cdots & 0 \\
c\frac{\alpha r_2}{r_1+\alpha N_r} & c & \cdots & 0 \\
\vdots & \vdots & \ddots & \vdots \\
c\frac{\alpha r_n}{r_1+\alpha N_r} & 0 & \cdots & c
\end{pmatrix},$$
$$D=\begin{pmatrix}\frac{cx_1(r_1+\alpha N_r-r_1)}{(r_1+\alpha N_r)^2}-b & -\frac{cx_1r_1\alpha}{(r_1+\alpha N_r)^2} & \cdots & -\frac{cx_1r_1\alpha}{(r_1+\alpha N_r)^2} \\
-\frac{c\alpha x_1r_2}{(r_1+\alpha N_r)^2} & \frac{c\alpha x_1(r_1+\alpha N_r-\alpha r_2)}{(r_1+\alpha N_r)^2}-b & \cdots & -\frac{c\alpha^2x_1r_2}{(r_1+\alpha N_r)^2} \\
\vdots & \vdots & \ddots & \vdots \\
-\frac{c\alpha x_1r_n}{(r_1+\alpha N_r)^2} & -\frac{c\alpha^2x_1r_n}{(r_1+\alpha N_r)^2} & \cdots & \frac{c\alpha x_1(r_1+\alpha N_r-\alpha r_n)}{(r_1+\alpha N_r)^2}-b\end{pmatrix}.$$

The fixed points, where the node 1 is persistent and all other nodes are neutral, are
$$x_i>0,1\le i\le n, r_1=0, r_i>0,2\le i\le n.$$
From these conditions we get
$$f_1=p\beta N_r=\beta\sum_{j=2}^nf_j, f_i=pr_i,0<x_1<\frac{b}{c}N_r,x_i=\frac{b}{c}r_i-\frac{r_i}{N_r}x_1,2\le i\le n.$$
So these fixed points belong to the subspace $f_1=\beta\sum_{j=2}^nf_j$.
At each fixed point, 
$$A_J=\mathbf{0},C=\begin{pmatrix}
0 & 0 & \cdots & 0 \\
\frac{cr_2}{N_r} & c & \cdots & 0 \\
\vdots & \vdots & \ddots & \vdots \\
\frac{cr_n}{N_r} & 0 & \cdots & c
\end{pmatrix},D=\begin{pmatrix}
\frac{cx_1}{\alpha N_r}-b & 0 & \cdots & 0 \\
-\frac{cx_1r_2}{\alpha N_r^2} & \frac{cx_1(N_r-r_2)}{N_r^2}-b & \cdots & -\frac{cx_1r_2}{N_r^2} \\
\vdots & \vdots & \ddots & \vdots \\
-\frac{cx_1r_n}{\alpha N_r^2} & -\frac{cx_1r_n}{N_r^2} & \cdots & \frac{cx_1(N_r-r_n)}{N_r^2}-b
\end{pmatrix}.$$
The matrix $B$ is an upper triangular matrix with negative diagonal entries.
Assume 
{now that}
the Jacobian is invertible. Then we get
$$\begin{pmatrix}A_J & B \\ C & D\end{pmatrix}\begin{pmatrix}
E & F \\ G & H \end{pmatrix}=\begin{pmatrix} 0 & B \\ C & D \end{pmatrix}\begin{pmatrix}
E & F \\ G & H\end{pmatrix}=\begin{pmatrix}
BG & BH \\ CE+DG & CF+DH \end{pmatrix}=\begin{pmatrix} I & 0 \\ 0 & I \end{pmatrix}.$$
Then $$BG=I,G=B^{-1};BH=0,H=0;CE+DG=0;CF+DH=CF=I.$$
Observe that $C$ is a lower triangular matrix with diagonal entries $0,c,c,\dots, c$.   This matrix is not invertible because there is no $F$ such that $CF=I$. 
Therefore we have shown that 
the Jacobian is not invertible. Hence it has a zero eigenvalue at
{any fixed point of the considered family.}

Now we compute the eigenvector which corresponds to the eigenvalue 0.
$$\begin{pmatrix} 0 & B \\ C & D \end{pmatrix}\begin{pmatrix} u \\ v \end{pmatrix}=\begin{pmatrix} Bv \\ Cu+Dv \end{pmatrix}=\begin{pmatrix} 0 \\ 0 \end{pmatrix}.$$
$$Bv=0\Rightarrow v=0; Cu+Dv=Cu=0, u=\begin{pmatrix} N_r & -r_2 & -r_3 & \cdots & -r_n \end{pmatrix}^T.$$

{Therefore the eigenspace corresponding to eigenvalue 0}
is such that all $r_i$'s are fixed (since $v=0$), but all the $x_i$'s can move along the the $u$ direction. The $x_i$ values at the fixed points are $0<x_1<\frac{b}{c}N_r,x_i=\frac{b}{c}r_i-\frac{x_1}{N_r}r_i,2\le i\le n$, and the $r_i$ values are $r_1=0, r_i=f_i/p,2\le i\le n$. So all the fixed points lie on the eigenspace of eigenvalue zero.

Clearly there is only one nonzero entry in the $(n+1)$-th row of the Jacobian at the fixed points. It is the $(n+1)$-th element $\frac{cx_1}{\alpha N_r}-b$ and it is an eigenvalue of the Jacobian. Thus this eigenvalue $\lambda_1=\frac{cx_1}{\alpha N_r}-b$ is negative when $x_1<\alpha\frac{b}{c}N_r$ and it is positive when $\alpha\frac{b}{c}N_r<x_1<\frac{b}{c}N_r$.

Therefore, as we have shown, neutral viruses can support the existence of persistent viruses. However, in contrast to altruistic viruses, neutral viruses can only support the existence of bounded concentrations of persistent viruses, while altruistic viruses can support any concentrations of persistent viruses \cite{BS}. Moreover, without the presence of altruistic viruses the states with persistent viruses could only be  marginally stable, while in the presence of altruistic viruses the state of local immunodeficiency can be stable \cite{BS}.

\section{Attaching minimal networks to large CR networks with random structure (topology)}

\begin{figure}[h]
    \vspace{-10pt}
    \includegraphics[scale=0.2]{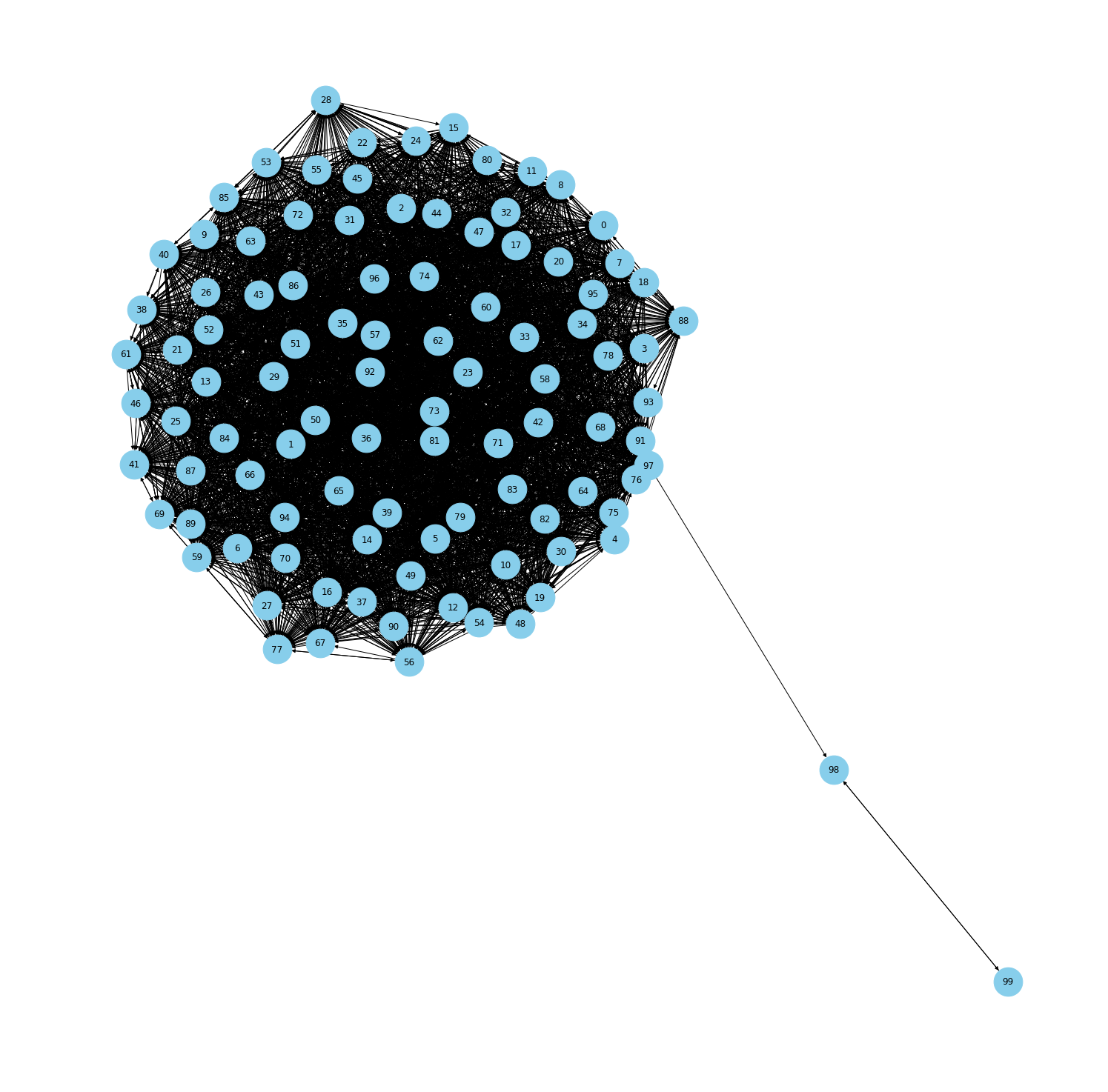}
    \caption{Attaching the minimal network to a random network}
    \label{dandelion}
\end{figure}

We
{have}
found two minimal CR networks with three
{types of viruses (nodes) each.}
In this section we analyze whether these networks can  maintain the
{state of}
local immunodeficiency
{in their persistent viruses}
if we attach them to a large random CR network.

{This question is important for understanding what may happen when CR networks are growing, particularly when some interactions between two hosts result in transmission of some CR subnetwork from one host to the other. In this section we discuss results of our numerical experiments attempting to address this question.}

{Shown in  Fig. \ref{dandelion} is a 100-node network that looks like a dandelion or a ball with a tail. The ball consists of nodes 0 through 97 and edges between them. The edges are generated by a $98\times98$ random matrix with 0-1 entries. The (i,j) entry of the matrix being 1 means there is an edge from node i-1 to node j-1 and 0 means no edge. The tail is built by identifying nodes 97, 98 and 99 with nodes 1, 2 and 3 from the branch-cycle CRN (Fig. \ref{minimal}), i.e. there is an edge from node 97 to 98, an edge from node 98 to 99 and an edge from node 99 back to node 98. Node 97 is where our minimal network, the tail, attaches to a large network, the ball. Since the edges inside of the ball are randomly generated, all nodes are equally important, attaching the minimal network to node 97 is as good as attaching it to any other node in the ball.}

{After building a randomly generated ball and attaching a minimal network to it as a tail, we build the evolution model \eqref{population} on this 100-node network. Then we solve it numerically by the forward Euler method. In each numerical experiment, we randomly generate a set of parameters. The parameters $f_i$'s are allowed to be different, and they are generated as a random vector with entry values between 0 and 1. The ranges for all parameters are listed below \eqref{parameters}. Initial conditions for virus and antibody populations are also randomly generated.}

\begin{equation}\label{parameters}
\begin{split}
f_i\sim U(0,1),f_{98}\sim U(1,2),p\sim U(0,1), b,c\sim U(0,5),\\
\alpha\sim U(0.5,1),\beta=\alpha^2,x_i(0),r_i(0)\sim U(0,0.1).
\end{split}
\end{equation}
Here $f_i\sim U(0,1)$ means that $f_i$ is a uniformly distributed random number in the interval $(0,1)$.
{The ranges of these parameters are consistent with the ones accepted in literature \cite{nowak00, rong, dahari, pnas, BS}.}

{We repeated this numerical experiment multiple times. The results were consistent in all simulations, even though the structure of the large network (the ball) changes every time (since it's randomly generated), and the parameters are different and randomly generated each time as well. We always observed local immunodeficiency at equilibrium on the tail (i.e. the attached minimal network).}

{We also repeated the same numerical experiments for the new minimal symmetric CRN (Fig. \ref{symmetric}). Again, we consistently observed local immunodeficiency at equilibrium on the tail. In conclusion, both minimal networks (one from \cite{BS} and the new one from the present paper) maintain local immunodeficiency after being attached to a large CR network.}

\section{Discussion}\label{dis}

\cite{BS} proved that local immunodeficiency discovered in \cite{pnas} is a stable and robust phenomenon which may already appear in CRNs with just three viruses. Therefore LI should likely be present in all diseases that demonstrate cross-immunoreactivity. Indeed, it is not necessary for CRNs to be large, which is typical for Hepatitis C \cite{pnas}, to have local immunodeficiency. We also rigorously demonstrated in \cite{BS} that there that there are easy ways to build larger networks with several persistent nodes (viruses) which remain invisible to the host's immune system because of their positions in the CRN. 
 
In the present paper we prove that the simplest cross-immunoreactivity network with three nodes (viruses) and just two edges can have a state of stable and robust local immunodeficiency. It is truly the smallest network of this type because no network with three nodes can have fewer than two edges. (It is worthwhile to mention that \cite{BS} showed no CRN with two nodes can maintain stable local immunodeficiency). These results raise a natural question whether it is possible to experimentally build such small  cross-immunoreactivity networks with stable and robust local immunodeficiency.

We also analyze here the role of neutral viruses in the creation of local immunodeficiency. It turns out that, in the absence of altruistic viruses, neutral viruses (even when any number of them is present to help a persistent node) can maintain local immunodeficiency only as a marginally stable state. Moreover, if only neutral and persistent viruses are present then the population of persistent viruses cannot exceed a relatively small value in a sharp contrast to the situation when altruistic viruses are also present and then population of persistent viruses becomes virtually unbounded \cite{pnas, BS}.

{Lastly we also showed numerically that small networks with stable LI can be attached to large CR networks and keep the state of stable local immunodeficiency.}

{Overall, it is shown that local immunodeficiency will likely be present in all diseases with cross-immunoreactivity. Indeed, we proved that very small CR networks may have a stable state with local immunodeficiency. Moreover, if such a small network is attached to another CR network this small sub-network may retain its stable LI state.}
All these results call for future numerical, analytic and, first of all, biological studies. The most important and pressing question is which types of viruses can play a role of persistent and/or altruistic ones \cite{NM}. Another question is whether one and the same virus may become persistent, altruistic or, perhaps, neutral depending on the structure (topology) of the corresponding cross-immunoreactivity network.


\appendix

\section{Computation of $P(\lambda)$ from the characteristic polynomial of the Jacobian for Fig. \ref{symmetric}}
$$P(\lambda)=-\lambda\begin{vmatrix}-\lambda & -p\beta x_3 & -px_3 \\ \frac{c\alpha r_2}{\alpha r_2+r_3} & \frac{b\alpha r_3}{\alpha r_2+r_3}-b-\lambda & -\frac{b\alpha r_2}{\alpha r_2+r_3} \\ \frac{cr_3}{\alpha r_2+r_3} & -\frac{b\alpha r_3}{\alpha r_2+r_3} & -\frac{br_3}{\alpha r_2+r_3}-\lambda\end{vmatrix}-p\beta x_1\begin{vmatrix} 0 & -\lambda & -px_3 \\ c & \frac{c\alpha r_2}{\alpha r_2+r_3} & -\frac{b\alpha r_2}{\alpha r_2+r_3} \\ 0 & \frac{cr_3}{\alpha r_2+r_3} & -\frac{br_3}{\alpha r_2+r_3}-\lambda\end{vmatrix}$$
$$=-\lambda\begin{vmatrix} -\lambda & -p\beta x_3 & -px_3 \\ \frac{c\alpha r_2}{\alpha r_2+r_3} & \frac{b\alpha r_3}{\alpha r_2+r_3}-b-\lambda & -\frac{b\alpha r_2}{\alpha r_2+r_3} \\ c & -b-\lambda & -b-\lambda\end{vmatrix}+cp\beta x_1\begin{vmatrix} -\lambda & -px_3 \\ \frac{cr_3}{\alpha r_2+r_3} & -\frac{br_3}{\alpha r_2+r_3}-\lambda\end{vmatrix}$$
$$=cp\beta x_1(\lambda^2+\frac{br_3}{\alpha r_2+r_3}\lambda+\frac{pcx_3r_3}{\alpha r_2+r_3})+\lambda^2\begin{vmatrix} \frac{b\alpha r_3}{\alpha r_2+r_3}-b-\lambda & -\frac{b\alpha r_2}{\alpha r_2+r_3} \\ -b-\lambda & -b-\lambda \end{vmatrix}$$
$$-\lambda p\beta x_3\begin{vmatrix} \frac{c\alpha r_2}{\alpha r_2+r_3} & -\frac{b\alpha r_2}{\alpha r_2+r_3} \\ c & -b-\lambda \end{vmatrix}+\lambda px_3\begin{vmatrix} \frac{c\alpha r_2}{\alpha r_2+r_3} & \frac{b\alpha r_3}{\alpha r_2+r_3}-b-\lambda \\ c & -b-\lambda \end{vmatrix}$$
$$=bf_1(1-\alpha)(\lambda^2+\frac{br_3}{\alpha r_2+r_3}\lambda+bpr_3)+\lambda^2(-b-\lambda)(\frac{b\alpha r_3}{\alpha r_3+r_3}-b-\lambda+\frac{b\alpha r_2}{\alpha r_2+r_3})$$$$-\lambda p\beta x_3(-\frac{c\alpha r_2}{\alpha r_2+r_3}\lambda-\frac{bc\alpha r_2}{\alpha r_2+r_3}+\frac{bc\alpha r_2}{\alpha r_2+r_3})+\lambda px_3(-\frac{c\alpha r_2}{\alpha r_2+r_3}\lambda-\frac{bc\alpha r_2}{\alpha r_2+r_3}-\frac{bc\alpha r_3}{\alpha r_2+r_3}+bc+c\lambda)$$
$$=bf_1(1-\alpha)(\lambda^2+\frac{bf_3}{\alpha r_2+r_3}\lambda+bpr_3)+\lambda^2(\lambda+b)(\lambda+b-\frac{b\alpha(r_2+r_3)}{\alpha r_2+r_3})+\lambda p\beta x_3\frac{c\alpha r_2}{\alpha r_2+r_3}\lambda$$
$$+px_3\lambda[\frac{cr_3}{\alpha r_2+r_3}\lambda+\frac{bcr_3(1-\alpha)}{\alpha r_2+r_3}]$$
$$=bf_1(1-\alpha)(\lambda^2+\frac{br_3}{\alpha r_2+r_3}\lambda+bpr_3)+\lambda^2(\lambda+b)(\lambda+\frac{br_3(1-\alpha)}{\alpha r_2+r_3})+\lambda^2p\beta\alpha r_2b+px_3\lambda[\frac{cr_3}{\alpha r_2+r_3}\lambda+\frac{bcr_3(1-\alpha)}{\alpha r_2+r_3}]$$
$$=\alpha bf_1\lambda^2+bf_1(1-\alpha)(\lambda^2+\frac{br_3}{\alpha r_2+r_3}\lambda+bpr_3)+\lambda^2(\lambda+b)(\lambda+\frac{br_3(1-\alpha)}{\alpha r_2+r_3})+px_3\lambda[\frac{cr_3}{\alpha r_2+r_3}\lambda+\frac{bcr_3(1-\alpha)}{\alpha r_2+r_3}].$$
\section{Computation of the characteristic polynomial of the Jacobian for Fig. \ref{2fig}}
Let $\lambda_1=\frac{cx_1}{\alpha r_2}-b$,
$$\det(J-\lambda I)=\det\begin{pmatrix}-\lambda & 0 & -px_1 & -p\beta x_1\\0 & -\lambda & 0 & -px_2\\ 0 & 0 & \lambda_1-\lambda & 0\\ c & c & -\frac{cx_1}{\alpha r_2} & -b-\lambda\end{pmatrix}=(\lambda_1-\lambda)\det\begin{pmatrix}-\lambda & 0 & -p\beta x_1 \\ 0 & -\lambda & -px_2 \\ c & c & -b-\lambda\end{pmatrix}$$
$$=(\lambda-\lambda_1)[\lambda(\lambda^2+b\lambda+cpx_2)+p\beta x_1c\lambda]=\lambda(\lambda-\lambda_1)[\lambda^2+b\lambda+cp(x_2+\beta x_1)]$$
$$=\lambda(\lambda-\lambda_1)P(\lambda),$$
\section{Computation of the characteristic polynomial of the Jacobian for Fig. \ref{3-node}}
Let $\lambda_1=\frac{cx_1}{\alpha(r_2+r_3)}-b<\frac{b}{\alpha}-b$, then
$$|J-\lambda I|=\begin{vmatrix} -\lambda & 0 & 0 & -px_1 & -p\beta x_1 & -p\beta x_1 \\
0 & -\lambda & 0 & 0 & -px_2 & 0 \\ 0 & 0 & -\lambda & 0 & 0 & -px_3 \\
0 & 0 & 0 & \lambda_1-\lambda & 0 & 0 \\
\frac{cr_2}{r_2+r_3} & c & 0 & -\frac{cx_1r_2}{\alpha(r_2+r_3)^2} & \frac{cx_1r_3}{(r_2+r_3)^2}-b-\lambda & -\frac{cx_1r_2}{(r_2+r_3)^2} \\
\frac{cr_3}{r_2+r_3} & 0 & c & -\frac{cx_1r_3}{\alpha(r_2+r_3)^2} & -\frac{cx_1r_3}{(r_2+r_3)^2} & \frac{cx_1r_2}{(r_2+r_3)^2}-b-\lambda\end{vmatrix}$$
$$\stackrel{(r_6\to r_5+r_6)}{\longeq}\begin{vmatrix}-\lambda & 0 & 0 & -px_1 & -p\beta x_1 & -p\beta x_1 \\
0 & -\lambda & 0 & 0 & -px_2 & 0 \\ 0 & 0 & -\lambda & 0 & 0 & -px_3 \\
0 & 0 & 0 & \lambda_1-\lambda & 0 & 0 \\
\frac{cr_2}{r_2+r_3} & c & 0 & -\frac{cx_1r_2}{\alpha(r_2+r_3)^2} & \frac{cx_1r_3}{(r_2+r_3)^2}-b-\lambda & -\frac{cx_1r_2}{(r_2+r_3)^2} \\
c & c & c & -\frac{cx_1}{\alpha(r_2+r_3)} & -b-\lambda & -b-\lambda\end{vmatrix}.$$
By expanding  along the fourth row we get
$$|J-\lambda I|=(\lambda_1-\lambda)\begin{vmatrix} -\lambda & 0 & 0 & -p\beta x_1 & -p\beta x_1 \\
0 & -\lambda & 0 & -px_2 & 0 \\ 0 & 0 & -\lambda & 0 & -px_3 \\
\frac{cr_2}{r_2+r_3} & c & 0 & \frac{cx_1r_3}{(r_2+r_3)^2}-b-\lambda & -\frac{cx_1r_2}{(r_2+r_3)^2} \\
c & c & c & -b-\lambda & -b-\lambda\end{vmatrix}.$$
Expanding now along the second row we obtain
$$|J-\lambda I|=(\lambda-\lambda_1)[\lambda D_1(\lambda)+px_2 D_2(\lambda)],$$
where
$$D_1(\lambda)=\begin{vmatrix} -\lambda & 0 & -p\beta x_1 & -p\beta x_1 \\
0 & -\lambda & 0 & -px_3 \\ \frac{cr_2}{r_2+r_3} & 0 & \frac{cx_1r_3}{(r_2+r_3)^2}-b-\lambda & -\frac{cx_1r_2}{(r_2+r_3)^2} \\
c & c & -b-\lambda & -b-\lambda\end{vmatrix}$$
$$=-\lambda\begin{vmatrix}-\lambda & -p\beta x_1 & -p\beta x_1\\
\frac{cr_2}{r_2+r_3} & \frac{cx_1r_3}{(r_2+r_3)^2}-b-\lambda & -\frac{cx_1r_2}{(r_2+r_3)^2}\\
c & -b-\lambda & -b-\lambda\end{vmatrix}-px_3\begin{vmatrix}-\lambda & 0 & -p\beta x_1\\
\frac{cr_2}{r_2+r_3} & 0 & \frac{cx_1r_3}{(r_2+r_3)^2}-b-\lambda\\
c & c & -b-\lambda\end{vmatrix}$$
$$=-\lambda\begin{vmatrix}-\lambda & 0 & -p\beta x_1\\
\frac{cr_2}{r_2+r_3} & \frac{cx_1}{r_2+r_3}-b-\lambda & -\frac{cx_1r_2}{(r_2+r_3)^2}\\
c & 0 & -b-\lambda\end{vmatrix}+cpx_3[\lambda^2+(b-\frac{cx_1r_3}{(r_2+r_3)^2})\lambda+cp\beta x_1\frac{r_2}{r_2+r_3}]$$
$$=\lambda(\lambda+b-\frac{cx_1}{r_2+r_3})[\lambda^2+b\lambda+cp\beta x_1]+cpx_3[\lambda^2+(b-\frac{cx_1r_3}{(r_2+r_3)^2})\lambda+cp\beta x_1\frac{r_2}{r_2+r_3}],$$
$$D_2(\lambda)=\begin{vmatrix}-\lambda & 0 & 0 & -p\beta x_1 \\ 0 & 0 & -\lambda & -px_3 \\ \frac{cr_2}{r_2+r_3} & c & 0 & -\frac{cx_1r_2}{(r_2+r_3)^2} \\ c & c & c & -b-\lambda\end{vmatrix}=c\lambda[\lambda^2+(b-\frac{cx_1r_2}{(r_2+r_3)^2})\lambda+cpx_3]$$
$$+cp\beta x_1\lambda\frac{cr_3}{r_2+r_3}.$$

\section*{Acknowledgements}
We are indebted to P. Skums for valuable discussions. This work was partially supported by the NSF grant CCF-BSF-1664836 and by the NIH grant 1R01EB025022.


\bibliographystyle{elsarticle-num}


\end{document}